\documentclass [12pt]{article}
\usepackage{amssymb,amsmath,amsthm}
\usepackage{color}

\renewcommand{\span}{span}
\DeclareMathOperator{\ind}{ind}

\usepackage[cp1250]{inputenc}
\usepackage[active]{srcltx}

\newtheorem{thm}{Theorem}[section]
\newtheorem{prop}[thm]{Proposition}
\newtheorem{lem}[thm]{Lemma}
 \newtheorem{rem}[thm]{Remark}
 \newtheorem{cor}[thm]{Corollary}
\newtheorem{defn}[thm]{Definition}
\newtheorem{ex}[thm]{Example}
\DeclareMathOperator{\supp}{supp}
\newcommand{\C}{\mathbb{C}}
\newcommand{\R}{\mathbb{R}}

\newcommand{\N}{\mathbb{N}}
\newcommand{\T}{\mathbb{T}}

\newcommand{\D}{\mathbb{D}}

\newcommand{\la}{\lambda}

\long\def\comment#1{{}}

 \title{Self-adjoint operators associated with Hankel moment matrices}
 \author{Christian Berg and Ryszard Szwarc}

\begin{document}

 \maketitle

\begin{abstract} 
In a paper from 2016 D. R. Yafaev initiated a study of closable Hankel forms associated with the moments $(m_n)$ of a positive measure with infinite support on the real line.
If $m_n=o(1)$ Yafaev characterized the closure of the form based on earlier work on quasi-Carleman operators. We give a new proof of the description of the closure based entirely on moment considerations. 

The main purpose of the present paper is a description of the self-adjoint Hankel  operators
associated with closed Hankel forms  in the Hilbert space of square summable sequences. 

We do this not only in the case $m_n=o(1)$ studied by Yafaev but  also in two other cases, where the Hankel form is closable, namely if the moment sequence is indeterminate or if the moment sequence is determinate with finite index of determinacy.
\end{abstract}



\section{Introduction} In the following let $\mathcal M^*$ denote the set of positive measures with moments of any order and infinite support on the real line $\R$. The moment sequence  of $\mu\in\mathcal M^*$ is denoted
\begin{equation}\label{eq:mom}
m_n=m_n(\mu)=\int_{-\infty}^\infty x^n\,d\mu(x),\quad n=0,1,\ldots,
\end{equation} 
and we let
\begin{equation}\label{eq:Han}
\mathcal H=\mathcal H_\mu=\left(m_{k+l}\right)_{k,l=0}^\infty=\{m_{k+l}\}
\end{equation}
denote the corresponding Hankel matrix.
The set $\C[x]$ of polynomials in one variable with complex coefficients is contained in the complex Hilbert space $L^2(\mu)$.  

In this paper we shall focus on the question if the matrix $\mathcal H_\mu$ leads to a self-adjoint operator $H=H_\mu$ (bounded or unbounded) in the Hilbert space $\ell^2$ of square summable complex sequences $g=(g_n)_{n\ge 0}$. This holds in three cases:
$\mu$ concentrated on the open interval $(-1,1)$, $\mu$ indeterminate and $\mu$ determinate with finite index of determinacy. 

The standard orthonormal basis in $\ell^2$ is denoted $e_k, k=0,1,\ldots$.

 Denoting by $\mathcal F$ the set of sequences $g\in\ell^2$ with only finitely many non-zero entries, we have a linear map $A$ from $\mathcal F$ to $\C[x]$
\begin{equation}\label{eq:A}
Ag(x)=\sum_{k=0}^\infty g_k x^k,\quad g\in\mathcal F, x\in\C,
\end{equation}
and a sesquilinear form $Q$ defined on $\mathcal F\times\mathcal F$
\begin{equation}\label{eq:quadr}
Q(g,h)=\int_{-\infty}^\infty Ag(x) \overline{Ah(x)}\,d\mu(x)=\sum_{k,l=0}^\infty m_{k+l}g_k\overline{h_l}, 
\end{equation}
called the {\it Hankel form associated with the  sequence} $(m_n)$.

If we consider $A$ as a densely defined operator from $\ell^2$ to $L^2(\mu)$ with domain $\mathcal F$, we write it $(A_\mu,\mathcal F)$. Similarly we write the form $Q$ as $(Q,\mathcal F)$, when we need to specify its domain $\mathcal F$.

Although we shall only be concerned with sesquilinear forms \eqref{eq:quadr}, where $(m_n)$ is a moment sequence given in \eqref{eq:mom}, it is of historical interest to mention that the right-hand expression in \eqref{eq:quadr} is a well-defined hermitian form  for any real sequence $(m_n)$,  and a famous theorem of Hamburger, see \cite{Ak}, more precisely asserts that  the form is positive definite, i.e., 
$$
\sum_{k,l=0}^\infty m_{k+l}g_k\overline{g_l}>0,\quad g=(g_n)_{n\ge 0}\in\mathcal F\setminus\{0\},
$$ 
if and only if $(m_n)$ is of the form \eqref{eq:mom} with $\mu\in\mathcal M^*$.

Furthermore, for a complex sequence $(m_n)$ the right-hand side of \eqref{eq:quadr} is bounded in the norm of $\ell^2$, i.e., there exists a constant
$K>0$ such that
$$
\left|\sum_{k,l=0}^\infty m_{k+l}g_k\overline{h_l}\right|\le K||g||\,||h||,\quad g,h\in\mathcal F,
$$
 if and only if 
$$
m_n=\int_{\T} z^nh(z)\,dz,\quad n\ge 0,
$$
for a bounded measurable function $h$ on the unit circle $\T$, $dz$ denoting normalized Haar measure. This is a fundamental result of Nehari, \cite{N}.

For Hankel forms associated with moment sequences we have the following classical result of Widom.

\begin{thm}\label{thm:Widom}{\rm Widom \cite[Theorem 3.1]{W}}. Let $(Q,\mathcal F)$ denote the Hankel form given by the moment sequence \eqref{eq:mom}.
Then $Q$ is bounded in the norm of $\ell^2$ if and only if $m_n=O(1/n)$. 
\end{thm}

Yafaev \cite{Y2},\cite{Y3} examined closability of the form \eqref{eq:quadr} and established that it is closable together with some quasi-analyticity assumption on $(m_n)$ if and only if $m_n=o(1)$, or equivalently that $\mu$ is concentrated on the open interval $(-1,1)$. Berg and Szwarc \cite{B:S3} pointed out that the form is also closable if the measure $\mu$ is indeterminate, or if it is determinate with finite index of determinacy. For the definition of the previous concepts see Section 2.

Note that for $g\in\ell^2$ and $z\in\D:=\{z\in\C : |z|<1\}$ we can define
\begin{equation}\label{eq:Az}
Ag(z):=\sum_{k=0}^\infty g_k z^k.
\end{equation}
Clearly $Ag$ is a holomorphic function in $\D$ belonging to the Hardy space $H^2(\D)$, and satisfying 
\begin{equation}\label{eq:|z|<1}
|Ag(z)|\le \frac{||g||}{\sqrt{1-|z|^2}},\quad z\in\D.
\end{equation}
In particular $Ag(x)$ is a continuous function for $x\in (-1,1)$. Note that $Ag$ has two different meanings in \eqref{eq:A} and \eqref{eq:Az}, but it should be clear from the context which definition is used. They agree of course for $z\in \D$.

Lemma 2.1 in \cite{Y2} states that  the form $(Q,\mathcal F)$ is closable if and only if $(A_\mu,\mathcal F)$ is a closable operator from $\ell^2$ to $L^2(\mu)$, and we add that in the affirmative case their closures have the same domain. 
This is a simple consequence of the special case of \eqref{eq:quadr}
$$
Q(g,g)=||Ag||_{L^2(\mu)}^2,\quad g\in\mathcal F.
$$
The adjoint operator of $(A_\mu,\mathcal F)$ is given in \cite[Lemma 2.2]{Y2} as an operator from $L^2(\mu)$ to $\ell^2$ with domain
\begin{equation}\label{eq:adj1}
D(A_\mu^*)=\left\{f\in L^2(\mu):\left(\int_{-\infty}^\infty f(x) x^n\,d\mu(x)\right)\in\ell^2\right\} 
\end{equation}
and given by
\begin{equation}\label{eq:adj2}
(A_\mu^*f)_n=\int_{-\infty}^\infty f(x) x^n\,d\mu(x),\quad n\ge 0. 
\end{equation}

For the case of $\mu$ being concentrated on the interval $(-1,1)$ Yafaev proved the following.

\begin{thm}\label{thm:thYaf} {\rm Yafaev \cite[Theorem 3.4]{Y2}}. Let $\mu\in\mathcal M^*$ be concentrated on $(-1,1)$.
Then $(Q,{\mathcal F})$ is a closable form and the extension 
\begin{equation}\label{eq:Vmu0}
Q_\mu(g,h)=\int_{-1}^1 Ag(x)\overline{Ah(x)}\,d\mu(x),\quad g,h\in V_{\mu},  
\end{equation}
where $\mathcal F\subset V_{\mu}\subset \ell^2$ is defined as
\begin{equation}\label{eq:Vmu}
V_\mu:=\{g\in\ell^2 : Ag\in L^2(\mu)\},
\end{equation}
is a closed form $(Q_\mu,V_\mu)$,  which is the closure of $(Q,\mathcal F)$.
\end{thm}  

While it is easy to see that the form $(Q,\mathcal F)$ is closable and the form $(Q_\mu,V_\mu)$ is closed, it requires considerable work to see that the latter form is the closure of the former. The proof of this in \cite[pp 296--298]{Y2} is a reduction to an application of \cite[Theorem 3.9]{Y1} dealing with Laplace transforms.

At the end of this section we give a new proof of Theorem~\ref{thm:thYaf} which is independent of the results in \cite{Y1}. The main idea is contained in Lemma~\ref{thm:fund}. In Sections 3-5 we consider the main classes of measures $\mu\in\mathcal M^*$ for which the form $(Q,\mathcal F)$ is closable:
\begin{itemize}
\item Section 3: $\mu$ is concentrated on $(-1,1)$.
\item Section 4: $\mu$ is indeterminate.
\item Section 5: $\mu$ is determinate with finite index of determinacy.
\end{itemize}

In each section we study the {\it Hankel operator associated with $\mu$}, i.e., the positive self-adjoint operator $H=H_\mu$ in $\ell^2$ which corresponds to the positive closed form defined as the closure of $(Q,\mathcal F)$. We characterize the domain $D(H_\mu)$ of the operator and show how it acts on vectors $\xi\in D(H_\mu)$. The main result of Section 3 is Theorem~\ref{thm:o(1/n)}.

In Section 2 we give the necessary background about closable and closed forms and their corresponding self-adjoint operators. We also discuss relevant concepts from the theory of moments used in the following.

The difficult part of Yafaev's proof is the inclusion $V_\mu\subseteq D(\overline{A_\mu})=D(A_\mu^{**})$, i.e.,
\begin{equation}\label{eq:Yaf}
\int_{-1}^1 Ag(x)\overline{f(x)}\,d\mu(x)=\sum_{k=0}^\infty g_k\left(\int_{-1}^1 \overline{f(x)} x^k\,d\mu(x)\right),\quad g\in V_\mu, f\in D(A_\mu^*),
\end{equation}
which as mentioned is based on \cite{Y1}.

We start by giving a proof of  \eqref{eq:Yaf} in the framework of moment sequences, in particular without the use of \cite{Y1}.

\begin{proof} That $(Q_\mu,V_\mu)$ given by \eqref{eq:Vmu0},\eqref{eq:Vmu} is closed means that $V_\mu$ is complete in the inner product
\begin{equation*}\label{eq:V}
(g,h)_{V_\mu} =(g,h)_{\ell^2} +\int\limits_{-1}^1 Ag(x)\overline{Ah(x)}\,d\mu(x),
\end{equation*}
cf. Section 2.
To see the completeness,
 assume $(g^{(n)})$ is a Cauchy sequence in $V_\mu$. Then
there exist $g\in \ell^2$ and $f\in L^2(\mu)$ so that
$$
\|g^{(n)}-g\|_{\ell^2} \to 0,\qquad \|Ag^{(n)}-f\|_{L^2(\mu)}\to 0.
$$
This implies by \eqref{eq:|z|<1} that
 $$
Ag^{(n)} (x) \to \sum_{k=0}^\infty g_kx^k=Ag(x) 
$$
 uniformly on compact sets in $(-1,1).$ 
Therefore
 $$
Ag(x)=f(x),\quad \mu-{\rm almost\  everywhere},
 $$
and hence $g\in V_\mu$ and $\|g^{(n)}- g\|_{V_{\mu}}\to 0.$ 

We next claim that $\mathcal{F}$  is dense in $(V_\mu, (\cdot,\cdot)_{V_\mu})$, and prove it  below. It depends on the following result.

\begin{lem}\label{thm:fund} Assume that $\mu\in\mathcal M^*$ is concentrated on $(-1,1)$. For $g\in \ell^2$ the series $Ag(x)=\sum_{k=0}^\infty g_k x^k$ is absolutely convergent in $L^2(\nu),$ where $d\nu(x)=(1-x^2)d\mu(x).$ In particular $Ag\in L^2(\nu)$.
\end{lem}

 Indeed,
$$
 \sum_{k=0}^\infty |g_k| \, \|x^k\|_{L^2(\nu)}
\le \|g\|_{\ell^2} \left (\sum_{k=0}^\infty \|x^k\|^2_{L^2(\nu)}\right )^{1/2}.
$$
On the other hand
$$
\sum_{k=0}^\infty \|x^k\|^2_{L^2(\nu)}=\sum_{k=0}^\infty \int\limits_{-1}^1 x^{2k}(1-x^2)\,d\mu(x) =m_0,
$$
since $m_{2k}-m_{2k+2}\ge 0$, proving the lemma.

\medskip
To see that  $\mathcal{F}$ is  dense in $V_\mu$, let 
 $g\in V_\mu$ be perpendicular to $\mathcal F$, and we shall prove that $g=0$. By assumption $g\perp e_k$ for any $k,$ i.e.,
$$
 0=(g,e_k)_{V_\mu}=g_k+\int\limits_{-1}^1 Ag(x) x^k\,d\mu(x), \quad k\ge 0.
 $$
 This implies
 \begin{equation*}\label{eq:jump2}
 \int\limits_{-1}^1 Ag(x) x^k\,(1-x^2)\,d\mu(x)=g_{k+2}-g_k, \quad k\ge 0.
\end{equation*} 
 Using Lemma~\ref{thm:fund} we can multiply  sidewise by $\overline{g_k}$ and sum up to obtain
  $$
0\le \int\limits_{-1}^1 |Ag(x)|^2(1-x^2)\,d\mu(x)
  =\sum_{k=0}^\infty g_{k+2}\overline{g_k}-\|g\|_{\ell^2}^2.
$$
  Therefore
   $$
\left |\sum_{k=0}^\infty g_{k+2}\overline{g_k}\right|\ge \|g\|_{\ell^2}^2,
 $$
and since
$$
\left |\sum_{k=0}^\infty g_{k+2}\overline{g_k}\right|\le
\left(\sum_{k=0}^\infty |g_{k+2}|^2\right)^{1/2}||g||_{\ell^2}\le ||g||_{\ell^2}^2,
$$
we get equality in the Cauchy-Schwarz inequality, which implies $(g_{k+2})_{k\ge 0}=\alpha (g_k)_{k\ge 0}$ for some $\alpha\in\C$. If $g\neq 0$ then  
$|\alpha|=1$ and hence $|g_{k+2}|=|g_k|$ for all $k\ge 0$, but this is not possible for $g\in\ell^2\setminus\{0\}$.

Let us now prove Equation \eqref{eq:Yaf}.
For $g\in V_\mu$ there exists $g^{(n)}\in \mathcal{F}$ so that
 $\|g^{(n)}-g\|_{V_{\mu}}\to 0$.

 Equivalently
$$
\|g^{(n)}-g\|_{\ell^2}\to 0,\quad \|Ag^{(n)}-Ag\|_{L^2\mu)}\to 0.
$$
 For $f\in L^2(\mu)$ such that $(f_k)\in\ell^2$, where 
$f_k:=\int\limits_{-1}^1 f(x)x^k\,d\mu(x)$, 
we have 
\begin{eqnarray*}
\int\limits_{-1}^1 Ag(x)\overline{f(x)}\,d\mu(x)&=&\lim_n\int\limits_{-1}^1 Ag^{(n)}(x)\overline{f(x)}\,d\mu(x)\\
 &=&
\lim_n\sum_{k=0}^\infty g^{(n)}_k \overline{f_k}=\sum_{k=0}^\infty g_k\overline{f_k}.
\end{eqnarray*}
\end{proof}

\section{Preliminaries}

In the following we give the basic definitions and results about closable and closed forms
in the special case of positive forms. We refer to \cite{Bi:S} and \cite{R:S} for details.
 
Let  $E$ denote a complex Hilbert space with inner product $(x,y)$ and let $V\subseteq E$ denote a dense subspace equipped with a sesquilinear form $Q:V\times V\to\C$, i.e., $Q$ is linear in the first variable and conjugate linear in the second variable. We assume $Q$ to be a positive hermitian form meaning that  $Q(x,y)=\overline{Q(y,x)}$ and $Q(x,x)\ge 0$ for all $x\in V$.   We denote the form by the symbol $(Q,V)$ and call $V=D(Q)$ the domain of $Q$.

If we have two positive hermitian forms $(Q_i,V_i),i=1,2$, we write $(Q_1,V_1)\subseteq (Q_2,V_2)$ if $V_1\subseteq V_2$ and the restriction of $Q_2$ to $V_1\times V_1$ is equal to $Q_1$. 

\begin{defn}\label{thm:closed} A positive hermitian form $(Q,V)$ on $E$ is a  closed form, if $V$ is complete (hence a Hilbert space) under the inner
product
$$
(x,y)_V:=(x,y)+Q(x,y), \quad x,y\in V. 
$$
\end{defn}

There is  a one-to-one correspondence between positive hermitian closed forms $(Q,V)$ on a Hilbert space $E$ and positive self-adjoint operators $(H,D(H))$ on $E$, where $D(H)$ denotes the domain of the operator $H$. 

Firstly, if $(H,D(H))$ is a positive self-adjoint operator on $E$, then $H^{1/2}$ is a positive self-adjoint operator with domain $V:=D(H^{1/2})$ and $D(H)\subseteq
 V\subseteq E$. Therefore $V$ is dense in $E$ and 
\begin{equation*}\label{eq:sa1}
Q(x,y):=(H^{1/2}x,H^{1/2}y),\quad x,y\in V,
\end{equation*}
is a positive hermitian closed form on $V$ satisfying
\begin{equation*}\label{eq:sa2}
Q(x,y)=(Hx,y),\quad x\in D(H), y\in V.
\end{equation*}

This follows from the spectral theorem.

Secondly:

\begin{thm}\label{thm:cfsa} Let $(Q,V)$ denote a positive hermitian closed form on a complex Hilbert space $E$. Then there is one and only one positive self-adjoint operator $(H,D(H))$ on $E$ with $D(H)\subseteq V$ such that
\begin{equation*}\label{eq:cf1}
Q(x,y)=(Hx,y),\quad x\in D(H), y\in V.
\end{equation*} 

We have $D(H^{1/2})=V$ and
\begin{equation*}\label{eq:cf1a}
Q(x,y)=(H^{1/2}x,H^{1/2}y),\quad x,y\in V.
\end{equation*} 
\end{thm}

\begin{defn}\label{thm:closable} A positive hermitian form $(Q,V)$ on $E$ is called
closable if the two equivalent conditions of the following theorem hold.
\end{defn}

\begin{thm} The following two conditions on the positive hermitian form $(Q,V)$ on $E$ are equivalent:

(i) There exists a closed positive hermitian form $(Q',V')$ on $E$ such that $(Q,V)\subseteq (Q',V')$.

(ii) For any sequence $(x_n)$ from $V$ satisfying $||x_n||\to 0$ and $Q(x_n-x_m,x_n-x_m)
\to 0$ for $n,m\to\infty$, we have $\lim_{n\to\infty} Q(x_n,x_n)=0$.

If the two conditions hold,  then the restriction  $(Q',\overline{V})$, where $\overline{V}$ is the closure of $V$ in  $(V',(\cdot,\cdot)_{V'})$, is the smallest closed positive hermitian form extending $(Q,V)$ and is called its closure and denoted $(\overline{Q}, \overline{V})$. 
\end{thm}

The domain of the closure $(\overline{Q}, \overline{V})$ can be described as
$$
\overline{V}=\{x\in E : \exists (x_n) \in V, \lim_{n\to\infty}||x_n-x||=0,\lim_{n,m\to\infty}
Q(x_n-x_m,x_n-x_m)=0 \},
$$
and for $x,y\in\overline{V}$ we have
$$
Q(x,y)=\lim_{n\to\infty} Q(x_n,y_n),
$$
and this value is independent of the approximating sequences $(x_n),(y_n)$  from $V$ chosen for $x,y$.

\medskip
We next summarize some results connected to the moment problem with \cite{Ak},\cite{S} as basic references.

Associated with the moments \eqref{eq:mom} we have
the orthonormal polynomials $(P_n)$, which are  uniquely determined by the conditions
\begin{equation*}\label{eq:orpo}
\int_{-\infty}^\infty P_n(x)P_m(x)\,d\mu(x)=\delta_{n,m},
\end{equation*}
when we assume that all $P_n$ have  positive leading coefficients.
We write
\begin{equation}\label{eq:Pnxn}
P_n(x)=\sum_{k=0}^n b_{k,n}x^k,\quad x^n=\sum_{k=0}^n c_{k,n}P_k(x),
\end{equation}
and consider the infinite matrices $\mathcal B=\{b_{k,n}\},\mathcal C=\{c_{k,n}\}$,
where we define $b_{k,n}=c_{k,n}=0$ if $k>n$. It is well known that the coefficients $b_{k,n}, c_{k,n}$ are real.
 Clearly $\mathcal B\mathcal C=\mathcal C\mathcal B=\mathcal I$, where $\mathcal I=\{\delta_{k,n}\}$ is the identity matrix. Note also that ${\mathcal H}_\mu={\mathcal C}^t{\mathcal C}$ as matrix products with ${\mathcal H}_\mu$ given by \eqref{eq:Han}. 

If the moment problem is indeterminate, there exists an infinite convex set $M$ of measures $\mu$ satisfying \eqref{eq:mom}, but all these measures lead to the same orthonormal polynomials. All measures $\mu\in M$ have unbounded support. Among these measures are the Nevanlinna extremal or in short the  N-extremal (in \cite{S} called von Neumann solutions), which are precisely the measures $\mu\in M$ for which $\C[x]$ is dense in $L^2(\mu)$ by a theorem of M. Riesz (\cite[Chapters 2.2, 2.3]{Ak}). Each  N-extremal measure is  discrete and supported by the zero set $\Lambda$ of a certain entire function of minimal exponential type, i.e., of the form
$$
\mu=\sum_{\la\in\Lambda} c_\la \delta_{\la}, \quad c_{\la}>0,
$$ 
cf. \cite[p.101]{Ak}, where $\delta_{\la}$ denotes the Dirac measure with mass one concentrated in the point $\la\in\R$.

By a theorem going back to Stieltjes in special cases, the following remarkable fact holds: If one mass is removed from  $\mu$, then the new measure becomes determinate, i.e., 
$$
\widetilde{\mu}:=\mu-c_{\la_0}\delta_{\la_0},\quad \la_0\in\Lambda
$$ 
is determinate. For a proof see \cite[Theorem 7]{B:C}. The measure $\widetilde{\mu}$ is a so-called {\it determinate measure with index of determinacy} 0. If further $n\ge 1$ masses are removed, we arrive at a {\it determinate measure $\mu'$ with index of determinacy $n$}. See Berg-Dur{\'a}n \cite[Theorem 3.6]{B:D}. The paper \cite{B:D} also contains an intrinsic characterization of such measures $\mu'$, which we describe now.

 For a determinate measure $\mu\in\mathcal M^*$ and $z\in\C$, the measure
$|t-z|^{2k}\mu$, which has the density $t\mapsto |t-z|^{2k}$ with respect to $\mu$, can be determinate or indeterminate when $k$ is a natural number. It is easy to see that if $|t-z|^{2k_0}\mu$ is indeterminate for some natural number $k_0$, then all the measures $|t-z|^{2k}\mu$ are indeterminate for natural numbers $k>k_0$.
Berg-Dur{\'a}n defined
$$
\ind_z(\mu)=\sup\{k\in \N_0 : |t-z|^{2k}\mu \;\mbox{is determinate}\}, 
$$
where $\N _0=\{0,1,2,\ldots\}$. If $\ind_{z_0}(\mu)=\infty$ for some $z_0\in\C$, then
$\ind_z(\mu)=\infty$ for all $z\in\C$ by Corollary 3.4 in \cite{B:D}, so all the measures $|t-z|^{2k}\mu$ are determinate for $z\in\C, k\in\N_0$. Also, if $\mu$ is non-discrete, then $\ind_z(\mu)=\infty$ for all $z\in \C$. In  \cite[Theorem 3.9]{B:D} it is furthermore proved that if $\ind_{z_0}(\mu)$ is finite for some $z_0\in\C$, then $\ind_z(\mu)$ is constant $k\ge 0$ for $z$ in the complement of the support of $\mu$ and constant equal to $k+1$ in the support of $\mu$. Moreover, in this case $\mu$ is obtained from an N-extremal measure by removing the mass at $k+1$ points of the support of this N-extremal measure. Such an N-extremal measure is highly non-unique by a perturbation result of Berg and Christensen, see \cite[Theorem  8]{B:C}. In \cite{B:D1} one defined $\ind(\mu)=\infty$ in the  first case and $\ind(\mu)=k$ in the second case. 

In the indeterminate case the series
\begin{equation*}\label{eq:P}
P(z):=\left(\sum_{n=0}^\infty |P_n(z)|^2\right)^{1/2},\quad z\in\C
\end{equation*}
converges uniformly on compact subsets of $\C$, cf. \cite[Theorem 1.3.2]{Ak}, so $P$ is a positive continuous function.

For any polynomial $p(z)=\sum g_k P_k(z), g=(g_k)\in\mathcal F$ and any compact subset
$K\subset\C$ we therefore have
\begin{equation}\label{eq:rep}
\sup_{z\in K} |p(z)|\le P_K ||g||_{\ell^2},\quad P_K:=\sup_{z\in K} P(z), 
\end{equation} 
and independent of the measure $\mu \in M$ we have 
$||p||_{L^2(\mu)}=||g||_{\ell^2}$.

In  the indeterminate case  the
polynomials $P_n$ form  an orthonormal basis in $L^2(\mu)$ for all the N-extremal solutions $\mu$, and for the other solutions $\mu$ they form an orthonormal basis in the closure $\overline{\C[x]}^{L^2(\mu)}$.

It follows that this closure is isometrically isomorphic as Hilbert space with  the space $\mathcal E$ of  functions of the form
\begin{equation}\label{eq:ent}
u(x)=\sum_{k=0}^\infty g_k P_k(x),\quad g\in\ell^2, \;x\in\R
\end{equation}  
under the norm 
$$
||u||_{L^2(\mu)}=||g||_{\ell^2},
$$
and the series above converges in $L^2(\mu)$ for any $\mu\in M$. However, because of \eqref{eq:rep} the series in \eqref{eq:ent} converges locally uniformly in $\C$ to an entire function, which is a representative of $u$.  By a theorem of M. Riesz,  cf. \cite[p.56]{Ak}, all the functions $u\in\mathcal E$ are entire of minimal exponential type. Furthermore, the inequality \eqref{eq:rep} holds with $p$ replaced by $u$.

Note that  if 
$u_n,u\in\mathcal E$ satisfies $u_n\to u$ in $L^2(\mu)$, then $u_n\to u$ uniformly on compact subsets of $\C$. By a classical result of complex analysis it  also holds
that $D^k u_n\to D^k u$ locally uniformly in $\C$ for any $k\in\N$.

For later use we note the following:
\begin{lem}\label{thm:nex} Let $\mu$ be an N-extremal solution to an indeterminate moment problem. For $f\in L^2(\mu)$ the entire function
\begin{equation}\label{eq:nex} 
u_f(z)=\sum_{k=0}^\infty (f,P_k)_{L^2(\mu)}P_k(z),\quad z\in\C
\end{equation}
is equal to $f$ on $\supp(\mu)$.

If $(f_n)_{n\ge 1}$ is a sequence from $L^2(\mu)$ with entire extensions $(u_{f_n})$ given by \eqref{eq:nex}, and if $f_n\to f$ in $L^2(\mu)$, then $u_{f_n}(z)\to u_f(z)$ locally uniformly for $z\in\C$.
\end{lem}

 \section{Hankel operators in case $m_n=o(1)$}

By \cite[Theorem 1]{Y2} the condition $m_n=o(1)$ is equivalent to $\mu$ being concentrated on the open interval $(-1,1)$.

We first treat a simple case of $m_n=o(1)$, where the Hankel matrix defines a bounded operator on $\ell^2$.

\begin{thm}\label{thm:ell1} For a moment sequence $(m_n)$ given by \eqref{eq:mom} the following conditions are equivalent:

(i) $||m||_1:=\sum_{k=0}^\infty |m_k|<\infty$.

(ii) $\sum_{k=0}^\infty m_{2k}<\infty$.
 
 (iii) $\mu$ is concentrated on $(-1,1)$ and $\int_{-1}^1(1-x^2)^{-1}d\mu(x)<\infty$.
 
If the equivalent conditions hold, then the Hankel matrix ${\mathcal H}_\mu=\{m_{k+l}\}$ defines a bounded operator $H_\mu$ in $\ell^2$ by
\begin{equation*}\label{eq:H}
(H_{\mu} g)_n=\sum_{k=0}^\infty m_{n+k}g_k,\quad  g\in\ell^2, n\ge 0. 
\end{equation*}
We have $||H_\mu||\le ||m||_1$ and $H_\mu$ is positive, self-adjoint and of trace class. 
\end{thm}

\begin{proof} It is clear that (i) implies (ii). Assume next that (ii) holds, and hence
$$
\int \sum_{k=0}^\infty x^{2k}\,d\mu(x)<\infty,
$$
but since the integrand is infinite for $|x|\ge 1$, we see that $\mu$ is concentrated on $(-1,1)$ and (iii) holds. Finally, if (iii) holds, so does (ii) and furthermore 
$$
|m_{2k+1}|\le \int_{-1}^1 |x^{2k+1}|\,d\mu(x)\le \int_{-1}^1 x^{2k}\,d\mu(x)=m_{2k},
$$
so (i) follows. 

Assume now that (i), (ii) and (iii) hold. Then the rows of the Hankel matrix belong to $\ell^1\subset\ell^2$, and therefore ${\mathcal H}_\mu g$ is a well-defined column for each $g\in\ell^2$. That each of these columns belong to $\ell^2$ is now a consequence of the Schur test, see e.g. \cite[E. 3.2.17]{GKP} because $\sum_{k=0}^\infty |m_{k+l}|\le ||m||_1$ for each $l\ge 0$. The trace of $H_\mu$ is the sum in (ii).
\end{proof}

We next describe the positive self-adjoint operator $(H_\mu,D(H_\mu))$ associated with the Hankel matrix \eqref{eq:Han} for an arbitrary moment sequence $m_n=o(1)$. We stress that this condition is equivalent to $\mu$ being concentrated on the open interval $(-1,1)$.

Let 
\begin{equation}\label{eq:F0}
\mathcal F_0:=\span\{v_k:=e_k-e_{k+2}, k\ge 0\},
\end{equation}
which is a subspace of $\mathcal F$ and $e_k\notin \mathcal F_0$ for all $k$. Note that the vectors $v_k, k\ge 0$ are linearly independent.

 It is easy to see that $\mathcal F_0$ is dense in $\ell^2$, for if $g\in\ell^2$ and $g\perp \mathcal F_0$,  then $g_k=g_{k+2}$ for all $k$ and this is only possible for $g\in\ell^2$ if $g=0$. 

It is useful to introduce the bounded shift operator $S$ by $(Sg)_k=g_{k+1}$ for $g\in \ell^2$. Then $Se_0=0, Se_k=e_{k-1}$ for $k\ge 1$ and $(I-S^2)$ maps $\ell^2$ injectively onto the dense subspace $(I-S^2)(\ell^2)$ of $\ell^2$. Furthermore,  $S^*e_k=e_{k+1}$ and $v_k=(I-S^{*2})e_k$.

The measure $\nu=(1-x^2)\mu$ belongs to $\mathcal M^*$ and the moments are
$m_n(\nu)=m_n-m_{n+2}$. Condition (iii) of Theorem 3.1 is clearly satisfied by $\nu$,
so by Theorem 3.1 the Hankel matrix $\mathcal H_\nu$ defines a bounded operator on $\ell^2$ denoted $H_\nu$.

\begin{thm}\label{thm:o(1/n)} Let $\mu\in\mathcal M^*$ be concentrated on $(-1,1)$ and let $d\nu(x)=(1-x^2)d\mu(x)$ with the associated bounded Hankel operator $H_\nu$ on $\ell^2$.

The operator $(H,\mathcal F_0)$ defined by linearity on $\mathcal F_0$ from \eqref{eq:F0} as
\begin{equation}\label{eq:H0}
Hv_k:=\left(m_{n+k}-m_{n+k+2}\right)_{n=0}^\infty
\end{equation}
is essentially self-adjoint. Its closure $\overline{H}$ is equal to the positive self-adjoint operator $H_\mu$ associated with the closed form $(Q_\mu,V_\mu)$ defined in Theorem~\ref{thm:thYaf}. 

The operator $\overline{H}$ can be described as
\begin{equation}\label{eq:Hmu}
D(\overline{H})=\{g\in\ell^2 : H_\nu g\in (I-S^2)(\ell^2)\}
\end{equation} 
and
\begin{equation}\label{eq:Hmu1}
\overline{H} g=(I-S^2)^{-1}(H_{\nu}g)=\sum_{l=0}^\infty S^{2l}H_{\nu}g,\quad g\in D(\overline{H}),
\end{equation}
where the series is  convergent in $\ell^2$.
\end{thm} 

\begin{rem}{\rm Summarizing we have for $g\in D(H_{\mu})$
$$
(H_{\mu}g)_n=\sum_{l=0}^\infty\left(\sum_{k=0}^\infty (m_{k+n+2l}-m_{k+n+2l+2})g_k\right),
$$
where the inner sum is absolutely convergent, and the outer sum is convergent but not necessarily absolutely convergent.
}
\end{rem}

\begin{proof}
 Note that $H_{\nu} e_k=H v_k$,  and $Hv_k$ defined in
\eqref{eq:H0} belongs to $\ell^1\subset \ell^2$. The operator $(H,\mathcal F_0)$ is clearly symmetric.

To see that the deficiency indices are zero, it is enough to prove that if $g\in\ell^2$ satisfies $g\perp (H-iI)v_k$ for all $k\ge 0$, then $g=0$.

By assumption
$$
0=((H-iI)v_k,g)_{\ell^2}=(H_{\nu} e_k,g)_{\ell^2}-i((I-S^2)^*e_k,g)_{\ell^2},\quad k\ge 0,
$$
which shows that
$$
(e_k, H_{\nu} g)_{\ell^2}+(e_k,i(I-S^2)g)_{\ell^2}=0, \quad k\ge 0,
$$
hence
$$
H_{\nu} g + i(I-S^2)g=0.
$$
Making the inner product with $g$ gives
$$
(g, H_\nu g)_{\ell^2}+i((I-S^2)g,g)_{\ell^2}=0,
$$ 
and in particular the imaginary part of this quantity is vanishing,  hence
$$
\Re((I-S^2)g,g)_{\ell^2}=0.
$$
In other words
$$
||g||_{\ell^2}^2=\Re(S^2g,g)_{\ell^2}\le |(S^2g,g)_{\ell^2}|\le ||S^2g||_{\ell^2}
 ||g||_{\ell^2}\le ||g||_{\ell^2}^2,
$$
so there is equality in the Cauchy-Schwarz inequality implying $S^2g=\alpha g$ for some $\alpha\in\C$ with $|\alpha|=1$ if $g\neq 0$. This leads to $|g_{k+2}|=|g_k|, k\ge 0$, contradicting 
$g\in\ell^2\setminus \{0\}$.

For $k\ge0, n\ge 1$ we define
\begin{equation*}\label{eq:gkn}
g^{(k,n)}=e_k-\frac{1}{n}\sum_{j=1}^n e_{k+2j}=\frac{1}{n}\sum_{j=1}^n (e_k-e_{k+2j})\in\mathcal F_0.
\end{equation*} 

Clearly 
\begin{equation}\label{eq:gkn1}
||g^{(k,n)}-e_k||_{\ell^2}^2=\frac{1}{n},
\end{equation}
and 
$$
Ag^{(k,n)}(x)=x^k\left(1-\frac{1}{n}\sum_{j=1}^n x^{2j}\right),
$$
hence by \eqref{eq:quadr}
\begin{eqnarray*}
\lefteqn{Q(g^{(k,n)},g^{(l,m)})-m_{k+l}}\\
&=&\int_{-1}^1 x^{k+l}\left[\left(1-\frac{1}{n}\sum_{j=1}^n x^{2j}\right)\left(1-\frac{1}{m}\sum_{j=1}^m x^{2j}\right)-1\right]\,d\mu(x).
\end{eqnarray*}
We can now estimate
\begin{eqnarray*}
\lefteqn{\left|Q(g^{(k,n)},g^{(l,m)})-m_{k+l}\right| }\\
&\le & \frac{1}{n}\int_{-1}^1\sum_{j=1}^n x^{2j}\left(1-\frac{1}{m}\sum_{j=1}^m x^{2j}\right)\,d\mu(x) + \frac{1}{m}\int_{-1}^1 \sum_{j=1}^m x^{2j}\,d\mu(x)\\
&\le & \frac{1}{n}\int_{-1}^1\sum_{j=1}^n x^{2j}\,d\mu(x)+\frac{1}{m}\int_{-1}^1\sum_{j=1}^m x^{2j}\,d\mu(x)\\
&=& \frac{1}{n}\sum_{j=1}^n m_{2j}+\frac{1}{m}\sum_{j=1}^m m_{2j},
\end{eqnarray*}
which tends to 0 for $n,m\to\infty$. Here we use the fact that if  a sequence $(\alpha_k)_{k\ge 0}$ of complex numbers tends to $\alpha$, then so does the sequence of averages.
 
Summarizing
$$
\lim_{n,m\to\infty} Q(g^{(k,n)},g^{(l,m)})=m_{k+l}=Q(e_k,e_l),
$$
which together with \eqref{eq:gkn1} shows that $(g^{(k,n)})_n$ is Cauchy in $V_\mu$ with the inner product $(g,h)_{\ell^2}+Q_\mu(g,h)$. Since $(Q_\mu,V_\mu)$ is closed by Theorem~\ref{thm:thYaf}, we know that $(g^{(k,n)})_n$ is convergent to an element in $V_\mu$, which by \eqref{eq:gkn1} must be $e_k$. This shows that 
$e_k\in\overline{{\mathcal F}_0}$ (closure in $V_\mu$), hence
$$
(Q, \mathcal F)\subseteq (\overline{Q_0}, \overline{{\mathcal F}_0})\subseteq (Q_\mu,V_{\mu}),
$$
where $Q_0$ is the restriction of $Q$ to $\mathcal F_0\times \mathcal F_0$,  and finally
$$
(\overline{Q}, \overline{\mathcal F})=(\overline{Q_0}, \overline{{\mathcal F}_0})=(Q_\mu,V_\mu),
$$
because  $(Q_\mu,V_\mu)$ is the closure of  $(Q,\mathcal F)$ by Theorem~\ref{thm:thYaf}.

We can now conclude that $H_\mu=\overline{H}$, because one can see that the self-adjoint operator $\overline{H}$ is associated with the closed form 
$(\overline{Q_0},\overline{\mathcal F_0})$, i.e.,
$$
(\overline{H}v,w)_{\ell^2}=\overline{Q_0}(v,w),\quad v\in D(\overline{H}), w\in\overline{\mathcal F_0},
$$ 
which follows from $(Hv,w)_{\ell^2}=Q_0(v,w)$ for $v,w\in\mathcal F_0$ by a limit procedure. We skip the details because
$H_\mu=\overline{H}$ also follows from the inclusion $H_\mu\subseteq\overline{H}$, which is established below.

Since $\overline{H}$ is self-adjoint, we know that $\overline{H}=H^*$, and therefore
$g\in D(\overline{H})$ if and only if there exists a sequence $a\in\ell^2$ such that
\begin{equation}\label{eq:a}
(g,Hv_l)_{\ell^2}=(a,v_l)_{\ell^2},\quad l\ge 0,
\end{equation}
and $\overline{H}g=a$ if \eqref{eq:a} holds.
Equivalently
$$
(g,H_\nu e_l)_{\ell^2}=(a,(I-S^{*2})e_l)_{\ell^2}=((I-S^2)a,e_l)_{\ell^2},\quad l\ge 0,
$$
i.e., $H_{\nu} g=(I-S^2)a$. This gives
$$
\sum_{l=0}^N S^{2l}H_{\nu} g =a-S^{2N+2}a,\quad \lim_{N\to\infty}||S^{2N+2}a||_{\ell^2}=0,
$$
which proves \eqref{eq:Hmu} and \eqref{eq:Hmu1}.

Let us now see that $H_\mu\subseteq\overline{H}$.

We first note that $L^2(\mu)\subset L^2(\nu)$ because for $f\in L^2(\mu)$
$$
\int_{-1}^1 |f(x)|^2\,d\nu(x)=\int_{-1}^1 |f(x)|^2(1-x^2)\,d\mu(x)\le \int_{-1}^1|f(x)|^2\,d\mu(x).
$$
If $g\in D(H_\mu)$ then $Ag\in L^2(\mu)$, and we have
\begin{eqnarray*}
(H_\mu g,v_k)_{\ell^2}=\int_{-1}^1 Ag(x)\overline{Av_k(x)}\,d\mu(x)=\int_{-1}^1 Ag(x)x^k\,d\nu(x)=(H_\nu g,e_k)_{\ell^2}.
\end{eqnarray*} 
This shows that
$$
((I-S^2)(H_\mu g))_k=(H_\mu g, v_k)_{\ell^2}=(H_\nu g)_k,\quad k\ge 0,
$$
hence $(I-S^2)(H_\mu g)=H_\nu g$,  so $g\in D(\overline{H})$ and $\overline{H}g=H_\mu g$, which shows the assertion. 
\end{proof}

\begin{rem}\label{thm:Q_0} {\rm 
Observe that in general if $(Q_i,V_i), i=1,2$ are closed forms  such that $(Q_1,V_1)\subseteq (Q_2,V_2)$, there is no simple relation between the corresponding self-adjoint operators $H_1$ and $H_2$. It is only if there is equality between the forms that a conclusion is possible, see \cite[p.279-280]{R:S}. 
}
\end{rem}

We next give another representation of the Hankel operator $H_\mu$.

\begin{thm}\label{thm:A*A} Let $\mu\in{\mathcal M}^*$ be concentrated on $(-1,1)$ and let $(A_\mu,V_\mu)$ denote the closed operator from $\ell^2$ to $L^2(\mu)$ defined by
\begin{equation*}\label{eq:Amu1}
A_{\mu}g(x)=Ag(x)=\sum_{k=0}^\infty g_kx^k,\quad V_\mu=\{g\in\ell^2 : Ag\in L^2(\mu)\}.
\end{equation*}
Then $H_{\mu}=A_{\mu}^* A_{\mu}$, i.e.,
\begin{equation}\label{eq:H1dom1}
D(H_\mu)=\left\{g\in V_\mu: \left(\int_{-1}^1 Ag(x) x^n\,d\mu(x)\right) \in\ell^2\right\}
\end{equation}
and
\begin{equation}\label{eq:H1val1}
(H_{\mu}g)_n=\int_{-1}^1 Ag(x)x^n\,d\mu(x), \quad g\in D(H_{\mu}), n\ge 0.
\end{equation}
\end{thm}

\begin{proof} Since $H_{\mu}$ is the positive self-adjoint operator in $\ell^2$ associated
with the closed form $(Q_\mu,V_\mu)$, we know that the following sequence belongs to $\ell^2$,
$$
(H_{\mu}g)_n=(H_\mu g, e_n)_{\ell^2}=\int_{-1}^1 Ag(x) x^n\,d\mu(x),\quad g\in D(H_\mu), n\ge 0
$$
where the second equality sign comes from \eqref{eq:Vmu0},  showing that $D(H_\mu)\subseteq D(A_\mu^*A_\mu)$ and
$H_{\mu}g=A_{\mu}^*A_{\mu}g$ for $g\in D(H_{\mu})$. Here it is important that $(A_\mu, V_\mu)$ is the closure of $A_\mu|\mathcal F$, so that these two operators have the same adjoint $A_\mu^*$  given by \eqref{eq:adj1} and \eqref{eq:adj2}.

Since $A_{\mu}^*A_{\mu}$ is a positive self-adjoint operator in $\ell^2$ extending $H_{\mu}$, they are equal. 

We can also see directly that $D(A_{\mu}^*A_{\mu})\subseteq D(H_{\mu})$ in the following way: If $g\in D(A_{\mu}^*A_{\mu})$, then
$$
u_n:=\int_{-1}^1 Ag(x) x^n\,d\mu(x),\quad n\ge 0
$$
defines a sequence $u=(u_n)\in\ell^2$. Therefore
$$
u_n-u_{n+2}=\int_{-1}^1 Ag(x) x^n\,d\nu(x),
$$    
where $d\nu(x)=(1-x^2)d\mu(x)$, and we recall by Lemma~\ref{thm:fund} that
the series $\sum g_kx^k$ converges absolutely to $Ag$ in $L^2(\nu)$.
This gives for $N\to\infty$
$$
\sum_{k=0}^N g_k m_{n+k}(\nu)=(\sum_{k=0}^N g_kx^k, x^n)_{L^2(\nu)}\to
(Ag, x^n)_{L^2(\nu)}=u_n - u_{n+2},
$$
so $(H_\nu g)_n=u_n-u_{n+2}$ for $n\ge 0$,
hence $(I-S^2)u=H_{\nu}g$, showing that $g\in D(H_\mu)$.
\end{proof}

\begin{rem}{\rm The transition matrix $\mathcal{C}=\{c_{k,l}\}$ defined in \eqref{eq:Pnxn} 
 is the matrix of the linear mapping $A:{\mathcal F}\to \C[x]$, cf. \eqref{eq:A}, with respect to the bases $(e_k)$ and $(P_k)$, where the latter are the orthonormal polynomials with respect to $\mu$.

 We can also consider $\mathcal{C}$ as a densely defined operator in $\ell^2$ with domain $\mathcal F$ given by $g\mapsto {\mathcal C}g$.  Note that ${\mathcal C}(\mathcal F) \subseteq \mathcal F$.

Letting $U:L^2(\mu)\to\ell^2$ denote the isometric isomorphism determined by $U(P_k)=e_k, k\ge 0$, we have
\begin{equation*}\label{eq:U}
U(Ag)={\mathcal C}g,\quad g\in\mathcal F,
\end{equation*}
hence
\begin{equation*}\label{eq:C}
\int_{-1}^1 Ag(x)P_k(x)\,d\mu(x)=({\mathcal C}g)_k, \;||Ag||_{L^2(\mu)}=||{\mathcal C}g||_{\ell^2},\quad g\in\mathcal F, k\ge 0.
\end{equation*}
We can now conclude that the operator $(\mathcal C,\mathcal F)$ has a closure $\overline{\mathcal C}$ because $(A_{\mu}, \mathcal F)$ is closable and
$$
D(\overline{\mathcal C})=D(\overline{A_{\mu}|\mathcal F})=V_{\mu}, \; \; (\overline{\mathcal C}g)_k=\int_{-1}^1 Ag(x)P_k(x)\,d\mu(x),\quad g\in V_{\mu}.
$$
 Moreover, ${\mathcal C}^*=(\overline{\mathcal C})^*$ exists and is densely defined. Finally ${\mathcal C}^*{\overline{\mathcal C}}$ is a self-adjoint positive operator in $\ell^2$.
}
\end{rem}

\begin{thm}\label{thm:C*C} Let $\mu\in{\mathcal M}^*$ be concentrated on $(-1,1)$ and let $(H_\mu,D(H_\mu))$ denote the associated Hankel operator from Theorem~\ref{thm:o(1/n)} and Theorem~\ref{thm:A*A}. Then
$H_{\mu}={\mathcal C}^*{\overline{\mathcal C}}$.
\end{thm}

\begin{proof}
In fact, for $g\in D(H_\mu)$ we have $Ag\in L^2(\mu)$ and by \eqref{eq:H1val1}
\begin{eqnarray*}
(H_{\mu}g)_n&=&\int_{-1}^1 Ag(x)x^n\,d\mu(x)=\sum_{k=0}^n c_{k,n}\int_{-1}^1 Ag(x)P_k(x)\,d\mu(x)\\
&=&\sum_{k=0}^n c_{k,n} (\overline{\mathcal C}g)_k,
\end{eqnarray*}
hence
$$
(\overline{\mathcal C}g, {\mathcal C}e_n)_{\ell^2}=(H_{\mu}g,e_n)_{\ell^2},\quad n\ge 0, 
$$
showing that $\overline{\mathcal C}g\in D({\mathcal C}^*)$ and that
$$
{\mathcal C}^*(\overline{\mathcal C}g)=H_{\mu}g,
$$
so $H_{\mu}\subseteq {\mathcal C}^*\overline{\mathcal C}$. Since both operators are self-adjoint, we have equality.
\end{proof}

In case $(m_n)\in\ell^2$ the self-adjoint Hankel operator $H_{\mu}$ of Theorem~\ref{thm:o(1/n)} can be described in a simpler way because
the Hankel matrix ${\mathcal H}_{\mu}=\{m_{k+l}\}$ has square summable rows and columns. Therefore the matrix product ${\mathcal H}_{\mu}g$ defines a column of complex numbers for each $g\in\ell^2$ and ${\mathcal H}_{\mu}g\in\ell^2$ for $g\in\mathcal F$. This leads to the operator
${\mathcal H}_{\mu}$ with maximal domain
\begin{equation*}\label{eq:opell2}
 D_{max}:=\{ g\in\ell^2 : {\mathcal H}_{\mu}g\in\ell^2\}.
\end{equation*}
 
\begin{thm}\label{thm:ell2} Assume $m_n=\int x^n\,d\mu(x)$ is square summable.
Then the operator $({\mathcal H}_{\mu}, D_{max})$  is equal to the self-adjoint operator $H_{\mu}$.  
\end{thm}

\begin{proof} The densely defined operator ${\mathcal H}_{\mu}|{\mathcal F}$  is clearly symmetric. We claim that the deficiency indices are zero. To see this assume that
$g\in\ell^2$ is perpendicular to $({\mathcal H}_{\mu}-iI)(\mathcal F)$, i.e.,
$$
({\mathcal H}_{\mu} e_n, g)_{\ell^2}=i(e_n, g)_{\ell^2},\quad n\ge 0,
$$   
hence
$$
(H_{\nu}e_n, g)_{\ell^2}=({\mathcal H}_{\mu}(e_n-e_{n+2}), g)_{\ell^2}=
i((I-S^2)^*e_n, g)_{\ell^2},
$$
where $d\nu(x)=(1-x^2)\,d\mu(x)$ and $H_\nu$ is the associated bounded self-adjoint operator in $\ell^2$ as in  Theorem~\ref{thm:o(1/n)}.
This leads to
$$
H_{\nu}g+i(I-S^2)g=0,
$$
which like in the proof of Theorem~\ref{thm:o(1/n)} implies 
$\Re ((I-S^2)g, g)_{\ell^2}=0$ and hence
$g=0$.

Since $(\mathcal H_\mu|\mathcal F)^*=\overline{\mathcal H_\mu|\mathcal F}$ is self-adjoint, we see that $g\in D(\overline{\mathcal H_\mu|\mathcal F})$ if and only if $g\in\ell^2$ and there  exists $a\in\ell^2$ such that
$$
(\mathcal H_\mu g)_n=\sum_{j=0}^\infty m_{n+j}g_j= (g,\mathcal H_\mu|\mathcal F e_n)_{\ell^2}=(a,e_n)_{\ell^2}=a_n,\quad n\ge 0
$$
or equivalently if $g\in D_{max}$. Hence $D_{max}=D(\overline{\mathcal H_\mu|\mathcal F})$ so $(\mathcal H_\mu,D_{max})$ is self-adjoint. If $g\in\mathcal F$, then
by \eqref{eq:Vmu0}
$$
(H_\mu g)_n=(H_\mu g, e_n)_{\ell^2}=\int_{-1}^1 Ag(x)x^n\,d\mu(x)=(\mathcal H_\mu g)_n,\quad n\ge 0,
$$
which yields $g\in D_{max}$. Thus $H_\mu|\mathcal F_0\subset H_\mu|\mathcal F\subset \mathcal H_\mu|D_{max}$. Since $H_\mu=\overline{H_\mu|\mathcal F_0}$
and $\mathcal H_\mu|D_{max}$ are self-adjoint, they are equal.
\end{proof}

\begin{ex}\label{thm:ex1} {\rm For $c>0$ let $\mu_c$ denote the probability measure on the open interval $(0,1)$ with density $(1/\Gamma(c))(\log(1/x))^{c-1}$ with respect to Lebesgue measure.
The moments are
\begin{equation*}\label{eq:c-mom}
m_n(c)=\frac{1}{\Gamma(c)}\int_0^1 x^n(\log(1/x))^{c-1}\,dx=\frac{1}{(n+1)^c},\quad n\ge 0.
\end{equation*}

The corresponding Hankel matrix ${\mathcal H}_{\mu_c}$ leads to a positive self-adjoint operator $H_{\mu_c}$ since $m_n(c)=o(1)$. It is unbounded for $0<c<1$ and bounded for $c\ge 1$ by Theorem~\ref{thm:Widom}. For $c=1$, ${\mathcal H}_{\mu_1}$ is the Hilbert matrix and $H_{\mu_1}$ is of norm $\pi$, see \cite{Ch} and references therein.

For $c>1$, $H_{\mu_c}$ is of trace class by Theorem~\ref{thm:ell1}. For $c>1/2$ the operator can be described as in Theorem~\ref{thm:ell2}.

Let us consider the subspaces $V_\mu$ and $D(H_\mu)$ given in \eqref{eq:Vmu} and 
 \eqref{eq:H1dom1} in case $\mu=\mu_c$ with $0<c<1$
$$
V_{\mu_c}=\left\{g\in\ell^2 :\sum_{k=0}^\infty g_kx^k \in L^2(\mu_c)\right\}
$$
and
$$
D(H_{\mu_c})=\left\{g\in V_{\mu_c} : \left(\int_0^1 (\sum_{k=0}^\infty g_kx^k)x^n\,d\mu_c(x)\right)\in\ell^2\right\}.
$$
Denoting $g_d:=(1/(n+1)^d)\in\ell^2$ for $d>1/2$, it is not difficult to prove that
\begin{eqnarray*}
g_d \in V_{\mu_c}  &\iff& d>1-\frac{c}{2}\\
g_d \in D(H_{\mu_c})  &\iff& d>\frac{3}{2}-c,
\end{eqnarray*} 
so for $d\in(1-c/2,3/2-c]$ we have $g_d\in V_{\mu_c}\setminus  D(H_{\mu_c})$.  
}
\end{ex}

\begin{ex}\label{thm:ex2} {\rm For $\la>-1/2$ we consider the probability density
$$
w_{\la}(x)=\frac{1}{B(1/2,\la+1/2)} (1-x^2)^{\la-1/2},\quad -1 < x <1.
$$
The corresponding orthogonal polynomials are the Gegenbauer polynomials, see \cite{A:A:R}. The odd moments vanish  and the even moments are
$$
m_{2k}(\la)=\int_{-1}^1 x^{2k} w_{\la}(x)\,dx=\frac{(1/2)_k}{(\la+1)_k},
$$
where $(a)_k=a(a+1)\cdots (a+k-1)$. Note that
$$
m_{2k}(\la) \sim \frac{\Gamma(\la+1)}{\sqrt{\pi}}\frac{1}{k^{\la+1/2}},\quad k\to\infty.
$$
For $\la>-1/2$ the corresponding Hankel matrix defines a positive self-adjoint operator $H_\la$ in $\ell^2$, unbounded for $-1/2<\la<1/2$, bounded for $\la\ge 1/2$ and of trace class  for $\la>1/2$.
}
\end{ex}

\section{The indeterminate case}
Suppose now that $m_n=\int x^n\,d\mu(x)$ is an indeterminate moment sequence, and
that $\mu$ is an arbitrary solution from the infinite set $M$ of measures satisfying \eqref{eq:mom}. 

 As described in \cite{B:S3}, the matrices $\mathcal{B},\mathcal{C}$ introduced in \eqref{eq:Pnxn}, define closable operators in $\ell^2$ with domain $\mathcal F$. The closure $\overline{\mathcal B}$ is a Hilbert-Schmidt operator and it maps $\ell^2$  injectively onto a dense subspace $V:=\overline{\mathcal B}(\ell^2)$ of $\ell^2$, 
cf. \cite[Prop. 4.3]{B:S1}. The closure $\overline{\mathcal C}$ is the inverse of $\overline{\mathcal B}$ with domain $D(\overline{\mathcal C})=V$.

The operator ${\mathcal B}^*$ is also Hilbert-Schmidt and injective with dense range ${\mathcal B}^*(\ell^2)$ and its inverse
$({\mathcal B}^*)^{-1}={\mathcal C}^*$. In \cite[Section 4]{B:S1} we have discussed the operator $\mathcal A:=\overline{\mathcal B}{\mathcal B}^*$, which is an injective self-adjoint  operator of trace class. If  we identify $\mathcal A$ with the matrix $({\mathcal A}e_l, e_k)_{\ell^2}$ which equals the matrix product ${\mathcal B}{\mathcal B}^t$, it was discussed in \cite{B:S2} if the matrix product ${\mathcal A}{\mathcal H}_\mu$ exists in the sense that the series defining its elements are absolutely convergent, and if the product is equal to $\mathcal I$. This is true for some indeterminate moment problems but not for all.

Note that $V$ only depends on $(m_n)$ and not on the particular solution $\mu$. For each $\mu\in M$ the operator $(A_\mu,\mathcal F)$ defined in \eqref{eq:A} is closable and the closure $\overline{A_\mu}$ has domain $D(\overline{A_\mu})=V$, cf. \cite[Theorem 3.2]{B:S3}. Indeed, for $\xi\in V$ we have $\xi=\overline{\mathcal B}g$ for a unique $g\in\ell^2$ and
\begin{equation}\label{eq:twoex}
\overline{A_\mu}\xi(x):=\sum_{k=0}^\infty g_kP_k(x)=\sum_{k=0}^\infty \xi_kx^k,\quad x\in\C,
\end{equation}  
belongs to $L^2(\mu)$, the first series converges in $L^2(\mu)$ and both series converge locally uniformly in $\C$ to a function belonging to $\mathcal E$ defined in \eqref{eq:ent}.
Since the functions in \eqref{eq:twoex} are independent of $\mu\in M$, we introduce the operator $(\mathfrak{A},V)$ from $\ell^2$ to $\mathcal E$ by 
 \begin{equation}\label{eq:frak}
\mathfrak{A}\xi(x)=\sum_{k=0}^\infty \xi_k x^k,\quad \xi\in V=\overline{\mathcal B}(\ell^2).
\end{equation}
Note that for $\xi\in\mathcal F \subset V$ we have $\mathfrak{A}\xi(x)=A\xi(x), \;x\in\C$. 
Furthermore, for all N-extremal solutions $\mu\in M$
$$
\{\mathfrak{A}\xi(x)1_{\supp(\mu)}(x) :\xi\in V\}=L^2(\mu).
$$

Similarly, for $\xi, \eta \in V$ the expression
\begin{equation}\label{eq:ind}
{\mathfrak Q}(\xi,\eta):=\int_{-\infty}^\infty {\mathfrak A}\xi(x) \overline{{\mathfrak A}\eta(x)}\,d\mu(x),\quad \mu\in M
\end{equation}
is independent of $\mu\in M$ and defines a positive hermitian form on $V\times V$. It is the closure of $(Q,\mathcal F)$ given in \eqref{eq:quadr} because of  \cite[Theorem 3.2]{B:S3}.

 We next describe the  positive self-adjoint Hankel operator $(H,D(H))$ corresponding to $({\mathfrak Q},V)$ in accordance with Theorem~\ref{thm:cfsa}. Note the analogy with Theorem~\ref{thm:C*C}.

\begin{thm}\label{thm:indet} Let $(m_n)$ be an indeterminate moment sequence and $\mu\in M$. The positive self-adjoint operator $H$ associated with the closed form $(\mathfrak{Q},V)$ is
 $H={\mathcal C}^*\overline{\mathcal C}$. We have
\begin{equation*}\label{eq:Hindet}
D(H)=D({\mathcal C}^*\overline{\mathcal C})=\left\{\xi\in V : \left(\int_{-\infty}^\infty \mathfrak{A}\xi(x) x^n\,d\mu(x)\right)\in\ell^2\right\}
\end{equation*}
\begin{equation*}\label{eq:Hindet1}
(H\xi)_n=\int_{-\infty}^\infty {\mathfrak A}\xi(x)x^n\,d\mu(x),\;n\ge 0, \quad \xi\in D(H).
\end{equation*}
Furthermore, $H={\mathcal A}^{-1}$.
\end{thm}

\begin{proof} Let $\xi\in D(H)$. Then $\xi\in V$ has the form $\xi=\overline{\mathcal B}g$ with $g=\overline{\mathcal C}\xi\in\ell^2$.
Note that by \eqref{eq:twoex} and \eqref{eq:frak}
$$
g_k=\int_{-\infty}^\infty {\mathfrak A}\xi(x) P_k(x)\,d\mu(x).
$$
Using this and \eqref{eq:ind} we find
 \begin{eqnarray*}
(\overline{\mathcal C}\xi,\mathcal{C}e_n)_{\ell^2}&=&\sum_{k=0}^n g_k c_{k,n}=
\sum_{k=0}^n c_{k,n} \int_{-\infty}^\infty {\mathfrak A}\xi(x) P_k(x)\,d\mu(x)\\
&=&\int_{-\infty}^\infty {\mathfrak A}\xi(x) x^n\,d\mu(x)={\mathfrak Q}(\xi,e_n)=(H\xi,e_n)_{\ell^2},
\end{eqnarray*}
showing that $\overline{\mathcal C}\xi\in D({\mathcal C}^*)$ and ${\mathcal C}^*(\overline{\mathcal C}\xi)=H\xi$, hence $H\subseteq {\mathcal C}^*\overline{\mathcal C}$, and there is finally equality, because both operators are self-adjoint.

We next have to show that
$$
\left\{\xi\in V : u_n:=\int_{-\infty}^\infty \mathfrak{A}\xi(x) x^n\,d\mu(x),\;(u_n)\in\ell^2\right\}\subseteq  D({\mathcal C}^*\overline{\mathcal C}).
$$
For $\xi$ belonging to the left-hand side, we have $\xi=\overline{\mathcal B}g$ with
$g=\overline{{\mathcal C}}\xi\in\ell^2$, hence
$$
g_k=\int_{-\infty}^\infty \mathfrak{A}\xi(x)P_k(x)\,d\mu(x)=\sum_{j=0}^k b_{j,k}u_j=
({\mathcal B}^*u)_k,
$$
showing that $g\in{\mathcal B}^*(\ell^2)=D({\mathcal C}^*)$, i.e., $\xi\in D({\mathcal C}^*\overline{\mathcal C}).$
Finally,
$$
{\mathcal A}^{-1}=(\overline{\mathcal B}{\mathcal B}^*)^{-1}=({\mathcal B}^*)^{-1}
(\overline{\mathcal B})^{-1}={\mathcal C}^*\overline{\mathcal C}=H.
$$
\end{proof}

\begin{rem}\label{thm:Hbelow} {\rm Note that $e_k\notin D(H), k\ge 0$. In fact, if we assume $e_k\in D(H)$ we get by Theorem~\ref{thm:indet} that
$$
(He_k)_n=\int_{-\infty}^\infty \mathfrak{A}e_k(x) x^n\,d\mu(x)=m_{k+n},
$$
showing that the $k$'th column in ${\mathcal H}_\mu$ is in $\ell^2$, but this is not possible in the indeterminate case as noticed in \cite[p. 3]{B:S2}.

Furthermore, 
\begin{equation}\label{eq:Hbelow}
(H\xi,\xi)_{\ell^2}\ge {\la}_{\infty}||\xi||_{\ell^2}^2,\quad \xi\in D(H),
\end{equation}
where ${\la}_\infty =\lim_{N\to\infty}{\la}_N$ is the limit of the smallest eigenvalue
${\la}_N$ of the truncated Hankel matrix $(m_{k+l})_{k,l=0}^N$. We recall that ${\la}_{\infty}>0$ characterizes the indeterminate case, see \cite{B:C:I}. For additional information about ${\la}_\infty$ see \cite{B:S1}.

To see \eqref{eq:Hbelow}, we note that for $\xi\in D(H)$
$$
(H\xi,\xi)_{\ell^2}=||\overline{\mathcal C}\xi||_{\ell^2}^2,
$$
and for $\xi\in\mathcal F$
$$
||{\mathcal C}\xi||_{\ell^2}^2=({\mathcal C}^t{\mathcal C}\xi, \xi)_{\ell^2}=\sum_{k,l=0}^\infty m_{k+l} \xi_k\overline{\xi_l}\ge {\la}_\infty ||\xi||_{\ell^2}^2,
$$
because $({\mathcal C}^t{\mathcal C}e_k, e_l)_{\ell^2}=m_{k+l}$.
We finally use that $\overline{\mathcal C}$ is the closure of $(\mathcal C,\mathcal F)$.
}
\end{rem}

\medskip

For $t\in\C$ the sequence  $(P_k(t))$ belongs to $\ell^2$ and as a special case of \eqref{eq:twoex}  we put $\xi(t)=(\xi_k(t))=\overline{\mathcal B}(P_k(t))\in V$,  hence
\begin{equation}\label{eq:twoex1}
{\mathfrak A}(\xi(t))(x)=\sum_{k=0}^\infty P_k(t)P_k(x)=\sum_{k=0}^\infty \xi_k(t) x^k,\quad x\in\C.
\end{equation}
It is an entire function of $(t,x)\in\C^2$.

\begin{cor}\label{thm:|t|<1} For $|t|<1$ we have $\xi(t)\in D(H)$ and $H(\xi(t))=(t^k)$.
\end{cor} 

\begin{proof} We have
\begin{eqnarray*}
( (P_k(t)), {\mathcal C} e_n)_{\ell^2}=\sum_{k=0}^n c_{k,n}P_k(t)=t^n=((t^k), e_n)_{\ell^2}, 
\end{eqnarray*}
and since $(t^k)\in\ell^2$ this shows that $(P_k(t))\in D({\mathcal C}^*)$ and that
${\mathcal C}^*(P_k(t))=(t^k)$. However,  since $(P_k(t))=\overline{\mathcal C}( \xi(t))$ and $H={\mathcal C}^*\overline{\mathcal C}$, we get the assertion.
\end{proof}

\begin{thm}\label{thm:core} Let $K$ be a subset of $\D$  such that $K$ contains an accumulation point of $K$. Then the subspace $\span\{\xi(t) : t\in K\}$ is a core for the operator $(H,D(H))$. 
\end{thm}

\begin{proof} The graph $G(H)=\{(\xi,H(\xi)) : \xi\in D(H)\}$ is a closed subspace of the product Hilbert space $\ell^2\times \ell^2$. The assertion of the theorem means that
$$
\Xi:=\span\{(\xi(t),H(\xi(t))) : t\in K\}
$$
is dense in $G(H)$. To see this assume that $(\eta,H(\eta))$ is orthogonal to $\Xi$ for some $\eta\in D(H)$, and let us prove that $\eta=0$. 

By assumption
$$
(\eta, \xi(t))_{\ell^2}+(H\eta,H(\xi(t)))_{\ell^2}=0,\quad t\in K,
$$
hence by Corollary~\ref{thm:|t|<1}
\begin{eqnarray*}
0 &=& (\eta,\overline{\mathcal B}(P_k(t)))_{\ell^2} + (H\eta,(t^k))_{\ell^2} =
({\mathcal B}^*\eta, (P_k(t)))_{\ell^2} +  (H\eta,(t^k))_{\ell^2}\\
&=& (\overline{\mathcal B}({\mathcal B}^*\eta) + H\eta,(t^k))_{\ell^2}=
(({\mathcal A} +H)\eta, (t^k))_{\ell^2},
\end{eqnarray*}
where $\mathcal A=\overline{\mathcal B}{\mathcal B}^*$ is the trace class operator mentioned in the beginning of this section.

In other words, the power series with coefficients given by the $\ell^2$-sequence $({\mathcal A}+H)\eta$ vanishes on the complex conjugate of the set $K$ and therefore
$({\mathcal A}+H)\eta=0$ and a fortiori 
$$
({\mathcal A}\eta,\eta)_{\ell^2} +(H\eta,\eta)_{\ell^2}=0.
$$
Since both operators are positive we conclude 
$$
0=({\mathcal A}\eta,\eta)_{\ell^2}=||{\mathcal B}^*\eta||_{\ell^2}^2,
$$
and since ${\mathcal B}^*$ is injective we get $\eta=0$.
\end{proof}

\section{Measures with finite index of determinacy}
As in the previous section $m_n=\int x^n\,d\mu(x)$ denotes an indeterminate moment sequence, and we assume now that $\mu$  is an arbitrary N-extremal solution. It is well-known that $\supp(\mu)$
is an infinite discrete set. Let $N\ge 1$ be a natural number and let $x_1<x_2<\cdots<x_N$ be points from the support of $\mu$. The measure
\begin{equation}\label{eq:detN}
\widetilde{\mu}=\mu-\sum_{j=1}^N \mu(\{x_j\})\delta_{x_j}
\end{equation}
is a determinate measure with index of determinacy $\ind(\widetilde{\mu})=N-1$, and any determinate measure with finite index of determinacy arises in this way as described
in Section 2.

We know that the operator $(A_{\widetilde{\mu}},\mathcal F)$ given by \eqref{eq:A}
is closable  by Theorem 1.4 in \cite{B:S3}. We shall now determine the closure.

\begin{thm}\label{thm:if1} Let $\widetilde{\mu}$ be determinate with $\ind({\widetilde{\mu}})=N-1$ and given by \eqref{eq:detN} for some $N\ge 1$. The closure of $(A_{\widetilde{\mu}},\mathcal F)$ is given by 
\begin{equation*}\label{eq:detN1}
D(\overline{A_{\widetilde{\mu}}})=V=\overline{\mathcal B}(\ell^2),
\end{equation*}
and for $\xi\in V$
\begin{equation}\label{eq:detN2}
\overline{A_{\widetilde{\mu}}}\xi(x)={\mathfrak A}\xi(x),\quad  x\in \supp(\widetilde{\mu}),
\end{equation}
where ${\mathfrak A}\xi$ is given by \eqref{eq:frak} corresponding to the indeterminate moment sequence $(m_n)$.

Furthermore, $\overline{A_{\widetilde{\mu}}}(V)=L^2(\widetilde{\mu})$.
\end{thm}

\begin{proof}
"$V\subseteq D(\overline{A_{\widetilde{\mu}}})$"

In fact, for $\xi\in V$, $\xi=\overline{\mathcal B}g, g\in\ell^2$, we have Equation \eqref{eq:twoex},
so for $n\in\N$ there exists $\xi^{(n)}\in\mathcal F$ so that  for $n\to\infty$
$$
\sum_{k=0}^n g_kP_k(x)=\sum_{k=0}^\infty \xi^{(n)}_k x^k\to {\mathfrak A}\xi(x)
$$
in $L^2(\mu)$ and a fortiori in $L^2(\widetilde{\mu})$.

Since $(g_0,g_1,\ldots,g_n,0,0,\ldots)\to (g_k)$ in $\ell^2$ and $\overline{\mathcal B}$ is a bounded operator, we have $\xi^{(n)}\to\xi$ in $\ell^2$, hence $\xi\in  D(\overline{A_{\widetilde{\mu}}})$ and \eqref{eq:detN2} holds.

\medskip
"$D(\overline{A_{\widetilde{\mu}}})\subseteq V$"

Assume $\xi\in D(\overline{A_{\widetilde{\mu}}})$. There exists $\xi^{(n)}\in\mathcal F$ such that $\xi^{(n)}\to\xi$ in $\ell^2$ and $A_{\widetilde{\mu}}\xi^{(n)}$ is convergent in 
$L^2(\widetilde{\mu})$. We shall prove that
\begin{equation}\label{eq:detN3}
\lim_{n\to\infty}\sum_{k=0}^\infty \xi^{(n)}_k x_j^k
\end{equation}
exists for $j=1,\ldots,N$, and this will imply convergence of $A_{\mu}\xi^{(n)}$ in $L^2(\mu)$, and consequently $\xi\in D(\overline{A_\mu})=V$. 

To see the convergence in \eqref{eq:detN3} we choose $N$  points $0<y_1<y_2<\cdots<y_N<1$, disjoint from $\supp(\widetilde{\mu})$. The measure
\begin{equation}\label{eq:modif}
\tau:=\widetilde{\mu}+\sum_{j=1}^N \mu(\{x_j\})\delta_{y_j}
\end{equation}
is indeterminate and N-extremal, but the moments of $\tau$ need not be equal to $(m_n)$. 
Let $A_\tau$ be the operator \eqref{eq:A} from $\mathcal F$ to $L^2(\tau)$ associated with $\tau$.
Using that 
$(1,y_j,y_j^2,\ldots)\in\ell^2$ we  see that
$$
\lim_{n\to\infty} \sum_{k=0}^\infty\xi^{(n)}_ky_j^k = \sum_{k=0}^\infty \xi_k y_j^k,\quad j=1,\ldots,N.
$$
Therefore $A_{\tau}\xi^{(n)}$ is convergent in $L^2(\tau)$, and since $\tau$ is N-extremal, this implies by Lemma~\ref{thm:nex} that $\sum_{k=0}^\infty \xi^{(n)}_k z^k$ converges locally  uniformly in $\C$  and in particular in the points $x_j, j=1,\ldots,N$, i.e., the convergence of \eqref{eq:detN3} holds.

That any $\widetilde{f}\in L^2(\widetilde{\mu})$ is of the form $\widetilde{f}=\overline{A_{\widetilde{\mu}}}\xi$ for some $\xi\in V$ is seen by extending $\widetilde{f}$ by 0 at the points $x_1,\ldots,x_N$ to a function $f\in L^2(\mu)$. Therefore
${\mathfrak A}\xi=f$ for some $\xi \in V$, hence
${\mathfrak A}\xi(x)=f(x)$ for $x\in\supp(\mu)$ and in particular for $x\in\supp(\widetilde{\mu})$, i.e.,
 $\widetilde{f}=\overline{A_{\widetilde{\mu}}}\xi$.
\end{proof}

The following Corollary of Theorem~\ref{thm:if1} is  also a Corollary of \cite[Theorem 3.2]{B:S3}.
\begin{cor}\label{thm:clQ} Let $\widetilde{\mu}$ be determinate with $\ind({\widetilde{\mu}})=N-1$ and given by \eqref{eq:detN} for some $N\ge 1$. Then the closure of the form $$
Q_{\widetilde{\mu}}(\xi,\eta)=\int_{-\infty}^\infty A\xi(x)\overline{A\eta(x)}\,d\widetilde{\mu}(x),\quad \xi,\eta\in\mathcal F
$$ 
 has domain 
$V=\overline{\mathcal B}(\ell^2)$ and
\begin{equation*}\label{eq:clQ}
\overline{Q_{\widetilde{\mu}}}(\xi,\eta)=\int_{-\infty}^\infty {\mathfrak A}\xi(x)
\overline{{\mathfrak A}\eta(x)}\,d{\widetilde{\mu}}(x),\quad \xi,\eta\in V.
\end{equation*}
\end{cor}

\begin{thm}\label{thm:multiple} If $\widetilde{\mu}$ is given by \eqref{eq:detN} under the assumption that
$x_j\in (-1,1), j=1,\ldots,N$, then the positive self-adjoint Hankel operators $H$ and $\widetilde{H}$ associated with $\mu$ and $\widetilde{\mu}$ have the same domain $D(H)$ and
\begin{equation}\label{eq:formH}
H=\widetilde{H}+\sum_{j=1}^N (1-x_j^2)^{-1}\mu(\{x_j\})P_{x_j},
\end{equation}
where $P_{x_j}$ is the orthogonal projection on the one-dimensional space spanned by the unit vector 
$$
v_{x_j}:=\sqrt{1-x_j^2}\sum_{k=0}^\infty x_j^k e_k.
$$

The kernel of the operator $\widetilde{H}$  equals $\span\{\xi(x_j),j=1,\ldots,N\}$, where
$\xi(x_j)$ is given by \eqref{eq:twoex1}.
\end{thm}

\begin{proof} The quadratic forms ${\mathfrak Q}$ and $\overline{Q_{\widetilde{\mu}}}$ satisfy for $\xi,\eta\in V$
\begin{eqnarray*}
{\mathfrak Q}(\xi,\eta)&=&\int_{-\infty}^\infty {\mathfrak A}\xi(x)\overline{{\mathfrak A}\eta(x)}\,d\mu(x)=\overline{Q_{\widetilde{\mu}}}(\xi,\eta)+\sum_{j=1}^N\mu(\{x_j\})
{\mathfrak A}\xi(x_j)\overline{{\mathfrak A}\eta(x_j)}\\
&=&\overline{Q_{\widetilde{\mu}}}(\xi,\eta)+\sum_{j=1}^N \mu(\{x_j\})\left(\sum_{k=0}^\infty\xi_kx_j^k\right)
\left(\sum_{k=0}^\infty\overline{\eta_k}x_j^k\right),
\end{eqnarray*}
so for $\xi\in D(H), \eta\in V$
$$
(H\xi,\eta)_{\ell^2}=\overline{Q_{\widetilde{\mu}}}(\xi,\eta)+\sum_{j=1}^N (1-x_j^2)^{-1}\mu(\{x_j\})(\xi,v_{x_j})_{\ell^2}(v_{x_j},\eta)_{\ell^2}.
$$
This together  with uniqueness in Theorem~\ref{thm:cfsa} shows \eqref{eq:formH}. 

For $j=1,\ldots,N$ we have
\begin{eqnarray*}
(\widetilde{H}(\xi(x_j)),\xi(x_j))_{\ell^2}&=&\int_{-\infty}^\infty |{\mathfrak A}(\xi(x_j))(x)|^2\,d\widetilde{\mu}(x)\\
&=&
\int_{-\infty}^\infty \left|\sum_{k=0}^\infty P_k(x_j)P_k(x)\right|^2\,d\widetilde{\mu}(x)=0,
\end{eqnarray*}
because the reproducing kernel $\sum_{k=0}^\infty P_k(x_j)P_k(x)$ vanishes on 
$\supp(\mu)\setminus\{x_j\}\supseteq\supp(\widetilde{\mu})$. (This is a consequence
of Parseval's equation for the orthogonal functions $1_{x_j}, 1_{x_k}$ in $L^2(\mu)$, when $x_j,x_k$ are two different points in the support of an N-extremal measure $\mu$.) This implies $\widetilde{H}(\xi(x_j))=0$, and therefore $\widetilde{H}$ vanishes on $\span\{\xi(x_j),j=1,\ldots,N\}$.

Assume next that $\widetilde{H}\xi=0$ for some $\xi\in D(\widetilde{H})=D(H)$. From \eqref{eq:formH} we get
$$
H\xi=\sum_{j=1}^N (1-x_j^2)^{-1}\mu(\{x_j\})(\xi,v_{x_j})_{\ell^2}v_{x_j},
$$
but from Corollary~\ref{thm:|t|<1} we have $\sqrt{1-x_j^2}H(\xi(x_j))=v_{x_j}$ and therefore
$$
H\xi=H\left(\sum_{j=1}^N (1-x_j^2)^{-1/2}\mu(\{x_j\})(\xi,v_{x_j})_{\ell^2}\xi(x_j)\right).
$$
We know that $H$ is injective by \eqref{eq:Hbelow} in Remark~\ref{thm:Hbelow}, and therefore
$$
\xi \in \span\{\xi(x_j),j=1,\ldots,N\}.
$$

\end{proof}

\begin{prop} Let the assumptions be as in Theorem~\ref{thm:multiple} and let
$$
W:=\{\xi\in D(H) : ({\mathfrak A}\xi)(x_j)=0, j=1,\ldots,N\}.
$$
Then $D(H)=W \oplus \span\{\xi(x_j), j=1,\ldots,N\}$, and for $\xi\in W$ we have
\begin{equation}\label{eq:ineq}
(\widetilde{H}\xi,\xi)_{\ell^2}=(H\xi,\xi)_{\ell^2}=\int_{-\infty}^\infty |{\mathfrak A}\xi(x)|^2\,d\widetilde{\mu}(x)=\int_{-\infty}^\infty |{\mathfrak A}\xi(x)|^2\,d\mu(x)\ge {\la}_{\infty} ||\xi||^2,
\end{equation}
where ${\la}_\infty>0$ is given in \eqref{eq:Hbelow}.
\end{prop}

\begin{proof} Note that for $i,j=1,\ldots,N$
$$
{\mathfrak A}(\xi(x_j))(x_i)=\sum_{k=0}^\infty P_k(x_j)P_k(x_i)=\mu(\{x_i\})^{-1}\delta_{i,j},
$$
showing that the vector space sum is direct. To see that the sum is all of $D(H)$ let $\xi\in D(H)$ and define $\alpha_j:=({\mathfrak A}\xi)(x_j)\mu(\{x_j\}), j=1,\ldots,N$. Then we have
$$
{\mathfrak A}\left(\xi-\sum_{j=1}^N\alpha_j \xi(x_j)\right)(x_i)=0,\quad i=1,\ldots,N,
$$
showing that $\xi-\sum_{j=1}^N\alpha_j \xi(x_j)\in W$.

Formula \eqref{eq:ineq} is now clear since $\mu$ and $\widetilde{\mu}$ only differs in the points $x_j$, where ${\mathfrak A}\xi(x)$ vanishes.
\end{proof}

\begin{rem}{\rm In the general case of $\widetilde{\mu}$ being given by \eqref{eq:detN} without any information about the location of the removed points, we consider the modified N-extremal measure $\tau$ defined in \eqref{eq:modif}. The domain of the positive self-adjoint Hankel operator $\widetilde{H}$ associated with $\widetilde{\mu}$ is equal to the domain of the positive self-adjoint Hankel operator $H=H_{\tau}$ associated with $\tau$ and
$$
H=\widetilde{H}+\sum_{j=1}^N (1-y_j^2)^{-1}\mu(\{x_j\})P_{y_j},
$$
i.e.,  $\widetilde{H}$ is a perturbation of $H$ with a self-adjoint operator of finite rank.

As before the kernel of $\widetilde{H}$ is of dimension $N$.
}
\end{rem}

{\bf Acknowledgment} The authors want to thank the referee for a careful reading of the manuscript and for valuable comments.

\noindent
Christian Berg\\
Department of Mathematical Sciences, University of Copenhagen\\
Universitetsparken 5, DK-2100 Copenhagen, Denmark\\
e-mail: {\tt{berg@math.ku.dk}}

\vspace{0.4cm}
\noindent
Ryszard Szwarc\\
Institute of Mathematics, University of Wroc{\l}aw\\
pl.\ Grunwaldzki 2/4, 50-384 Wroc{\l}aw, Poland\\ 
e-mail: {\tt{szwarc2@gmail.com}}

\end{document}